\documentclass[fleqn]{mat01}
\usepackage{times,mathtimy,amssymb,latexsym,amsbsy}
\begin{document}

\setcounter{page}{111}
\firstpage{111}

\newtheorem{theor}{\bf Theorem}

\def\coro{\trivlist\item[\hskip\labelsep{{\rm COROLLARY}}]}
\def\definit{\trivlist\item[\hskip\labelsep{{\rm DEFINITION}}]}
\def\remar{\trivlist\item[\hskip\labelsep{{\it Remarks.}}]}

\title{Measure free martingales}

\markboth{Rajeeva L~Karandikar and M~G~Nadkarni}{Measure free martingales}

\author{RAJEEVA L~KARANDIKAR and M~G~NADKARNI$^{*}$}

\address{Indian Statistical Institute, New Delhi 110 016, India\\
\noindent $^{*}$Department of Mathematics, University of Mumbai, Kalina,
Mumbai~400~098, India\\
\noindent E-mail: rlk@isid.ac.in; nadkarni@math.mu.ac.in}

\volume{115}

\mon{February}

\parts{1}

\pubyear{2005}

\Date{MS received 19 October 2004}

\begin{abstract}
We give a necessary and sufficient condition on a sequence of functions
on a set $\Omega$ under which there is a measure on $\Omega$ which
renders the given sequence of functions a martingale. Further such a
measure is unique if we impose a natural maximum entropy condition on
the conditional probabilities.
\end{abstract}


\maketitle

\section{Introduction}

The notion of measure free martingale is implicit in the
construction of equivalent martingale measures in the theory of
asset pricing in financial mathematics \cite{1,2}, but it has not
been fully isolated and made free of probability. Rather it has
remained hidden by specific processes and terminology of asset
pricing theory. We define a martingale purely in terms of sets and
functions, called measure free martingale, and show that every
martingale is a measure free martingale and conversely that every
measure free martingale admits a probability measure, which may be
finitely additive, under which it is a martingale. We describe the
convex set (together with their extreme points) of all probability
measures under which a measure free martingale is a martingale.
Among these measures there is one which in some sense is most
symmetric or most well spread, and entirely determined by the
measure free martingale. Boltzmann's entropy maximizing
distribution is needed here. To the best of our knowledge
probabilist's have not asked the simple question as to when a
sequence of function is a martingale under some measure. The
answer is relatively easy but has some pedagogic as well as
research value.

\section{Means of finite set of points}

Let $x_1, x_2, x_3,\ldots, x_k$ be $k$ real numbers, with repetitions
allowed. Assume that $x_1$ and $x_k$ are respectively the smallest and
the largest of $x_1, x_2,\ldots, x_k$. Let $\alpha$ be a real number.
Then there exists a probability vector $(p_1, p_2, \ldots, p_k)$ such
that
\begin{equation*}
x_1p_1 + x_2p_2 + \cdots + x_kp_k = \alpha,
\end{equation*}
if and only if $x_1 \leq \alpha \leq x_k$. If $k=2$ and $x_{1}
\neq x_{2}$, such a probability vector is unique. If $k>2$, it is
not unique without some additional requirements.

A result of Boltzmann proved using Lagrange's multipliers says
that there is a unique probability vector $(p_1, p_2,\ldots, p_k)$
which satisfies $x_1p_1 + x_2p_2 + \cdots + x_kp_k = \alpha$, and
maximizes the entropy
\begin{equation*}
-p_1{\log {p_1}} - p_2{\log{p_2}} -\cdots - p_k{\log{p_k}}.
\end{equation*}
It is given by
\begin{equation*}
p_j = \frac{\exp{(\lambda x_j)}}{\sum_{i=1}^k{\exp{(\lambda
x_i)}}}, \quad i = 1,2, \ldots, k,
\end{equation*}
where $\lambda$ is a constant.

We will call these probabilities the Boltzmann probabilities for
$x_1,x_2,\ldots, x_k; \alpha$.

In this connection it should be noted that for a fixed $x_1, x_2, \ldots,
x_k$ and variable $\lambda$, the probabilities
\begin{equation*}
p_i(\lambda) = \frac{\exp (\lambda x_i)}{\sum_{i=1}^k\exp (\lambda
x_i)},\quad i = 1, 2, \ldots, k
\end{equation*}
of $x_1, x_2, \ldots, x_k$ respectively  have the mean
$\sum_{i=1}^kx_ip_i(\lambda)$ which we denote by $m(\lambda)$.
Since $x_1$ and $x_k$ are minimum and maximum of $x_1, x_2,
\ldots, x_k$, we have
\begin{equation*}
\lim_{\lambda \rightarrow -\infty}p_i(\lambda) =
\frac{\delta_{1,i}}{n_i}, \quad \lim_{\lambda \rightarrow
\infty}p_{i}(\lambda) = \frac{\delta_{k,i}}{n_i},
\end{equation*}
where $n_i$ is the frequency of occurrence of $x_i$ in $x_1, x_2, \ldots,
x_k$. As a consequence,
\begin{equation*}
\lim_{\lambda \rightarrow -\infty }m(\lambda) = x_1,\quad \lim_{\lambda
\rightarrow \infty }m(\lambda) = x_k.
\end{equation*}
A calculation shows that ${\rm d}m/{\rm d}\lambda = v(\lambda)>
0$, where $v(\lambda)$ is the variance of the system $x_1, x_2,
\ldots, x_k$ with probabilities $p_1(\lambda), p_2(\lambda),
\ldots, p_k(\lambda)$. Thus $m(\lambda)$ is a strictly increasing
function of $\lambda$ which assumes every value between $x_1$ and
$x_k$. If $m(\lambda) = \alpha$, then $p_1(\lambda), p_2(\lambda),
\ldots, p_k(\lambda)$ are the probabilities which maximize the
entropy for the constraint $\sum_{i=1}^kp_k x_k = \alpha$. (See
\cite{3}, p.~172 for a related discussion of Boltzmann
distribution in the continuous case.)

Suppose $x_{1}, x_{2}, \ldots, x_{k}$ are distinct. The set $C$ of
probability vectors $(p_1, p_2, \ldots, p_k)$ such that
$\sum_{j=1}^kx_jp_j = \alpha$ is a convex set. It is easy to see
that its extreme points are precisely those $(p_1, p_2, \ldots,
p_k) \in C$ which have at most two non-zero entries.

\section{Measure free martingales}

Let $\Omega$ be a non-empty set. Let $f_n$, $n = 1,2,3,\ldots$ be a
sequence of real valued functions such that each $f_n$ has a finite range,
say $(x_{n1}, x_{n2}, \ldots, x_{nk_n})$, and these values are assumed
on the subsets $\Omega_{n1}, \Omega_{n2}, \ldots, \Omega_{nk_n}$. These
sets form a partition of $\Omega$ which we denote by ${\mathbb{P}}_n$.
We denote by ${\mathbb {Q}}_n$ the partition generated by ${\mathbb
{P}}_1, {\mathbb {P}}_2, \ldots, {\mathbb {P}}_n$ and the algebra
generated by ${\mathbb {Q}}_n$ is denoted by ${\cal{A}}_n$. Let
${\cal{A}}_\infty$ denote the algebra $\cup_{n=1}^\infty{\cal{A}}_n$.

Define ${\cal{A}}_n$ measurable functions $m_n,M_n$ as follows: For
$Q\in {\mathbb {Q}}_n$ and $\omega\in Q$,
\begin{align*}
m_n(\omega) &= \min_{q\in Q}f_{n+1}(q),\\[.2pc]
M_n(\omega) &= \max_{q\in Q}f_{n+1}(q).
\end{align*}

\vspace{.1pc}
\begin{definit}$\left.\right.$\vspace{.5pc}

\noindent The sequence $(f_n, {\cal{A}}_n)_{n=1}^\infty$ is said to
be a measure free martingale or probability free martingale if
\begin{equation*}
m_n(\omega)\leq f_n(\omega) \leq M_n(\omega) \quad \forall
\omega\in\Omega,\;n\geq 1.
\end{equation*}

Clearly, for each $Q\in {\mathbb{Q}}_n$, the function $f_n$ is
constant on $Q$. We denote this constant by $f_n(Q)$. With this
notation, it is easy to see that $(f_n, {\cal{A}}_n)_{n=1}^\infty$
is a measure free martingale or probability free martingale if and
only if for each $n$ and for each $ Q\in {\mathbb{Q}}_n$, $f_n(Q)$
lies between the minimum and the maximum values of
$f_{n+1}(Q^\prime)$ as $Q^\prime$ runs over $Q\cap
{\mathbb{Q}}_{n+1}$.

It is easy to see that if there is a probability measure on
${\cal{A}}_\infty$ with respect to which $(f_n,
{\cal{A}}_n)_{n=1}^\infty$ is a martingale, then $(f_n,
{\cal{A}}_n)_{n=1}^\infty$ is also a measure free martingale.
Indeed, let $P$ be such a measure. Then, for any $Q$ in
${\mathbb{Q}}_{n}$, $f_{n}(Q)$ is equal to
\begin{equation*}
\frac{1}{P(Q)}\sum_{ \{Q^\prime \in {\mathbb{Q}}_{n+1},
Q^\prime\subseteq Q\}} f_{n+1}(Q^\prime)P(Q^\prime),
\end{equation*}
so that $f_{n}(Q)$ lies between the minimum and the maximum values
$f_{n+1}(Q^\prime)$, $Q^\prime \in Q\cap {\mathbb{Q}}_{n+1}$. The
theorem below proves the converse.
\end{definit}

\begin{theor}[\!]
Given a measure free martingale $(f_n,
{\cal{A}}_n)_{n=1}^\infty${\rm ,} there exists for each $n\geq
0${\rm ,} a measure $P_n$ on ${\cal{A}}_n$ such that
\begin{equation*}
P_{n + 1} |_{{\cal{A}}_{n}} = P_{n}, \quad E_{n+1}(f_{n+1} |
{\cal{A}}_{n}) = f_{n},
\end{equation*}
where $E_{n+1}$ denotes the conditional expectation with respect to the
probability measure $P_{n+1}$. There is a finitely additive probability
measure $P$ on the algebra ${\cal{A}}_\infty${\rm ,} which may be countably
additive{\rm ,} such that for each $n${\rm ,} $P |_{{\cal{A}}_n} = P_n$.
\end{theor}

\begin{proof}
Define $P_1$ on ${{\cal{A}}_1}$ arbitrarily. Having defined $P_1, P_2,
\ldots, P_n $ on ${\cal{A}}_1,{\cal{A}}_2, \ldots, {\cal{A}}_n$ such
that
\begin{equation*}
P_j | {\cal{A}}_{j-1} = P_{j-1}, \quad E_j(f_j | {\cal{A}}_{j-1})
= f_{j-1}, \quad j = 2, 3, \ldots, n,
\end{equation*}
we define $P_{n+1}$ on ${\cal{A}}_{n+1}$ as follows: Choose an element
$Q$ in ${\mathbb{Q}}_{n}$. Let $A_1, A_2, \ldots, A_l$ be the partition
of $Q$ induced by $f_{n+1}$ so that $f_{n+1}$ assumes $l$ distinct
values, say $a_1, a_2,\ldots, a_l$, on $A_1, A_2, \ldots, A_l$
respectively. Let $a=f_n(Q)$ (the value assumed by $f_n$ on $Q$). Since
$(f_n, {\cal{A}}_n)_{n=1}^\infty$ is a measure free martingale, $a$ lies
between the minimum and the maximum values of $f_{n+1}$ on $Q$, so there
is a probability vector $(p_1, p_2, \ldots, p_l)$ such that
\begin{equation*}
a_1p_1 + a_2p_2+ \cdots + a_lp_l = a.
\end{equation*}

We define
\begin{equation*}
P_{n+1}(Q_i) = p_iP_n(Q),\quad i = 1, 2, \ldots, l.
\end{equation*}
Carrying out this procedure for all $Q \in {\mathbb{Q}}_n$ we get a
probability measure $P_{n+1}$ on ${\cal{A}}_{n+1}$ for which it is easy
to check that
\begin{equation*}
P_{n+1} | {\cal{A}}_n = P_n, \quad E_{n+1}(f_{n+1} | {\cal{A}}_n)
= f_n.
\end{equation*}

Induction completes the proof of the existence of the measures
$P_n$. Define $P$ by setting, for $A \in {\cal{A}}_\infty$, $P(A)
= P_n(A)$, if $A\in {\cal{A}}_n$. Thus the theorem stands
proved.\hfill $\Box$
\end{proof}

\begin{remar}
The measure $P$ on ${\cal{A}}_\infty$ may be called a martingale
measure associated to the measure free martingale $(f_n,
{\cal{A}}_n)_{n=1}^\infty$. The totality of such measures forms a
convex set whose extreme points are precisely those $P$ which have
the property that for any $n$ and for any $Q\in {\mathbb{Q}}_n$,
$P$ (hence $P_{n+1}$) assigns positive probability to at most two
elements in the partition of $Q$ induced by $f_{n+1}$. If, for
each $n$ and for each $Q\in {\mathbb{Q}}_n$, $Q\cap
{\mathbb{Q}}_{n+1}$ has two or less elements, then there is only
one martingale measure for the measure free martingale $(f_n,
{\cal{A}}_n)_{n=1}^\infty$.

Let $Q$ be an element in ${\mathbb{Q}}_n$. If we assign Boltzmann
probabilities of the values of $f_{n+1}$ on $Q$ to the
corresponding elements of the partition of $Q$ induced by
$f_{n+1}$, then we have the following theorem.
\end{remar}

\begin{theor}[\!]
Let $(f_n, {\cal{A}}_n)_{n=1}^\infty$ be a measure free
martingale. Then there is a unique probability measure $P$ on
${\cal{A}}_\infty$ such that\vspace{-.5pc}
\begin{enumerate}
\renewcommand\labelenumi{\rm (\arabic{enumi})}
\item $(f_n, {\cal{A}}_n)_{n=1}^\infty$ is a martingale with respect
to $P$.

\item For each $n$ and for each $Q\in {\mathbb{Q}}_n$ if $Q_1, Q_2,
\ldots, Q_l$ are the elements of $Q\cap {\mathbb{Q}}_{n+1}${\rm ,} then
$P(Q_1)/P(Q), P(Q_2)/P(Q), \ldots, P(Q_l)/P(Q)$ are the unique
probabilities which maximize
\begin{equation*}
\hskip -1.25pc -\sum_{i=1}^lp_i{\log {p_i}},
\end{equation*}
subject to the condition $\sum_{i=1}^la_ip_i = a${\rm ,} where $a$
is the value of $f_n$ on $Q$ and $a_1,a_2,\ldots, a_l$ are the
values assumed by $f_{n+1}$ on $Q$.

\item The probabilities $P(Q_i), i = 1,2,\ldots,l$ are given by the
formula{\rm :}
\begin{equation*}
\hskip -1.25pc P(Q_i) = P(Q)\cdot \frac{\exp({\lambda
a_i})}{\sum_{i=1}^l\exp({\lambda a_i})},
\end{equation*}
where $\lambda$ is a constant depending on $a, a_1, a_2, \ldots,
a_l$.\vspace{-.5pc}
\end{enumerate}
\end{theor}

In a certain sense this distribution $P$ of Theorem~2 may be viewed as
most symmetric or most well spread for the given measure free
martingale. It is determined entirely by the measure free martingale.
One may call $P$ the Boltzmann measure associated to the measure free
martingale $(f_n, {\cal{A}}_n)_{n=1}^\infty$, and the resulting measure
theoretic martingale, the Boltzmann martingale.

In the theory of asset pricing in financial mathematics there is
an important point of existence of equivalent martingale. Here, as
a consequence of Theorem~2, we have the following:\vspace{.7pc}

\begin{coro}$\left.\right.$\vspace{.5pc}

\noindent {\it With the notation of Theorem~$2$ above{\rm ,} if
$m$ is a probability measure on ${\cal{A}}_\infty$~ for which
there exist two positive constants $C$ and $D$ such that for all
$A\in \cup_{n=1}^\infty{\mathbb{Q}}_n${\rm ,}
\begin{equation*}
C \leq m(A)/P(A) \leq D,
\end{equation*}
then there is measure on ${\cal{A}}_\infty${\rm ,} {\rm (}e.g.{\rm ,}
$P${\rm ),} which is equivalent to $m$ and with respect to which $(f_n,
{\cal{A}}_n)_{n=1}^\infty$ is a martingale. This martingale measure is
unique{\rm ,} and equal to $P${\rm ,} if we require{\rm ,} for each $n$
and for each $Q \in {\mathbb{Q}}_n${\rm ,} the conditional distribution
on $Q\cap{\mathbb{Q}}_{n+1}$ to have maximum entropy.}\vspace{.7pc}
\end{coro}

A question arises. Note that we can associate the number $\lambda$
to the set $Q$ in Theorem~2. When we do this for all $Q\in
{\mathbb{Q}}_n$ we have a function $g_n$ defined on $\Omega$. Is
$g_n, n= 1,2,3, \ldots$ a measure free martingale?

Suppose $\Omega$ is a compact metric space and that sets in
${\cal{A}}_\infty$ form a clopen base for its topology. Then any
martingale measure for the measure free martingale $(f_n,
{\cal{A}}_n)_{n=1}^\infty$ extends to a countably additive measure on
the Borel field ${\cal{B}}$ of $\Omega$. The collection $C$ of all
martingale measures for $(f_n, {\cal{A}}_n)_{n=1}^\infty$ defined on
${\cal{B}}$ forms a compact convex set under weak topology, whose
extreme points are already described above.

\section{A result on convergence}

Let $\Omega$ be a compact metric space and let
$(f_n)_{n=1}^\infty$ be a sequence of continuous real valued
functions on $\Omega$. Let ${\mathbb {Q}}_n$ be the partition of
$\Omega$ generated by $f_1, f_2, \ldots, f_n$. Elements of
${\mathbb{Q}}_n$ are closed sets. Say that $(f_n,
{\mathbb{Q}}_n)_{n=1}^\infty$ is a martingale of continuous
functions if for each $n$ and for each $C\in {\mathbb{Q}}_n$ the
value of $f_n$ on $C$ lies between the minimum and the maximum
value of $f_{n+1}$ on $C$. We have the following theorem.

\begin{theor}[\!]
If the martingale $(f_n, {\mathbb{Q}}_n)_{n=1}^\infty$ of continuous
functions is also an equicontinuous sequence{\rm ,} i.e.{\rm ,} the
sequence of functions $(f_n)_{n=1}^\infty$ is equicontinuous{\rm ,} then
$(f_n)_{n=1}^\infty$ converges pointwise.
\end{theor}

\begin{proof}
Let ${\mathbb{Q}}_\infty$ denote the common refinement of all the
${\mathbb{Q}}_n, n = 1,2,\ldots$ and assume that
${\mathbb{Q}}_\infty$ is made of singleton sets. Let $\omega$ be a
point of $\Omega$ and let $C_n$ be the element of ${\mathbb{Q}}_n$
to which $\omega$ belongs. Then $\cap_{n=1}^\infty C_n =
\{\omega\}$, and since $C_n$'s are closed, we see that the
diameter of $C_n$ tends to zero as $n$ tends to $ \infty$. By
martingale and equicontinuity property of the sequence
$(f_n)_{n=1}^\infty$ we conclude that given any $\epsilon >0$
there is an $n_0$ such that for $n \geq n_0$, $ | f_n(\omega) -
f_{n_0}(\omega) | < \epsilon$. So $(f_n)_{n=1}^\infty$ converges
pointwise.

If ${\mathbb{Q}}_\infty$ is not made of singletons, then we
consider $\overline{\Omega} = \Omega/{\mathbb{Q}}_\infty$ equipped
with the quotient topology. Define for $c \in
{\mathbb{Q}}_\infty$, ${\overline{f}}_n(c) =$ the constant value
of $f_n$ on $c$. We can view ${\mathbb{Q}}_n$ also as a partition
of $\overline{\Omega}$. The sequence $( {\overline{f}}_n,
{\mathbb{Q}}_n)_{n=1}^\infty$ forms a martingale of continuous
functions on the compact set ${\overline{\Omega}}$ and the
functions ${\overline{f}}_n, n = 1,2,\ldots$ form an
equicontinuous sequence of functions. The common refinement
${\mathbb{Q}}_\infty$ of the partitions ${\mathbb{Q}}_n, n = 1, 2,
\ldots$ when considered as partition of $\overline{\Omega}$ is
the partition of $\overline{\Omega}$ into singleton sets. By
considerations of the previous paragraph we see that the sequence
$({\overline{f}}_n)_{n=1}^\infty$ converges pointwise, whence the
sequence $(f_n)_{n=1}^\infty$ converges pointwise. The theorem is
proved.\hfill $\Box$
\end{proof}

We conclude by raising a question about Boltzmann distribution.
Let $C$ be a compact subset of the real line and let $\alpha$ be
strictly between maximum and minimum points of $C$. Let $x_1, x_2,
\ldots, x_{k_\epsilon}$ be an $\epsilon$-net in $C$. Let
$\mu_\epsilon$ denote the Boltzmann distribution on this
$\epsilon$-net and $\alpha$. Can one say that $\mu_\epsilon$
converges weakly to a unique probability measure on $C$ as
$\epsilon \rightarrow 0$, independent of the choice of the
$\epsilon$-nets?

\section*{Acknowledgements}

It is a pleasure to thank R~B~Bapat for useful discussions. The
second author would like to acknowledge the hospitality provided
by IMSc, CMI and IFMR, Chennai during the period this note was
written.

\end{document}